\documentclass[12pt,twoside,reqno]{amsart}
\usepackage{amsfonts}
\usepackage{amssymb}
\usepackage[notcite,notref]{%showkeys
}

\makeatletter
\@namedef{subjclassname@2010}{%
  \textup{2010} Mathematics Subject Classification}
\makeatother

\numberwithin{equation}{section}

\newtheorem{teo}{Theorem}[section]

\theoremstyle{definition}

\newtheorem{rem}[teo]{Remark}

\numberwithin{equation}{section}

\def\R{\mathbb{R}}
\def\N{\mathbb{N}}

\def\e{\varepsilon}

\def\f{\varphi}

\begin{document}

%%%%%%%%%%%%%%%%%%%%%%%%%%%%%%%%%%%%%%%%%%%%%%%%%%%%%%%%%%%%%%%%%%%%%%

\title[The Hardy operator and its dual]{Optimal relations  between $L^P$-norms for the Hardy operator and its dual}

\author[V.I. Kolyada]{V.I. Kolyada}
\address{Department of Mathematics\\
Karlstad University\\
Universitetsgatan 1 \\
651 88 Karlstad\\
SWEDEN} \email{viktor.kolyada@kau.se}


\begin{abstract}
We obtain  sharp two-sided inequalities  between  $L^p-$\newline norms
$(1<p<\infty)$ of  functions $Hf$ and $H^*f$, where $H$ is the Hardy operator,  $H^*$ is its dual,
and $f$ is a nonnegative measurable function on $(0,\infty).$  In  an equivalent
form,  it gives sharp constants in the two-sided relations
 between $L^p$-norms
of functions $H\f-\f$ and $\f$, where $\f$ is a nonnegative
nonincreasing function on $(0,+\infty)$ with $\f(+\infty)=0.$ In
particular, it provides
 an alternative proof  of a result
obtained by N. Kruglyak and E. Setterqvist (2008) for $p=2k
\,\,(k\in \N)$ and by S. Boza and J. Soria (2011) for all $p\ge 2$,
and gives a sharp version of this result  for $1<p<2$.

\end{abstract}

%%%%%%%%%%%%%%%%%%%%%%%%%%%%%%%%%%%%%%%%%%%%%%%%%%%%%%%%%%%%%%%%%%%%%%

\maketitle

\date{}

\maketitle

\section{Introduction and main results}

Denote by $\mathcal M^+(\R_+)$ the class of all
 nonnegative measurable functions on $\R_+\equiv
(0,+\infty).$ Let  $f\in \mathcal M^+(\R_+).$ Set
$$
Hf(x)=\frac1x\int_0^x f(t)\,dt
$$
and
$$
H^*f(x)=\int_x^\infty \frac{f(t)}{t}\,dt.
$$
These equalities define the classical Hardy operator $H$ and its
dual operator $H^*$.  By Hardy's inequalities
\cite[Ch. 9]{HLP}, these operators are bounded in $L^p(\R_+)$ for any
$1<p<\infty.$ Furthermore, it is easy to show that for any $f\in \mathcal M^+(\R_+)$ and any $1<p<\infty$ the $L^p-$norms of
$Hf$ and $H^*f$ are equivalent. Indeed, let $f\in \mathcal M^+(\R_+).$
 By Fubini's theorem,
$$
Hf(x)=\frac1x\int_0^x dt\int_t^x\frac{f(u)}{u}\,du\le \frac1x\int_0^x H^*f(t)\,dt.
$$
On the other hand, Fubini's theorem gives that
$$
H^*f(x)=\int_x^\infty\frac{du}{u^2}\int_x^u f(t)\,dt\le \int_x^\infty\frac{Hf(u)}{u}\,du.
$$
Using these estimates and applying Hardy's inequalities \cite[p. 240, 244]{HLP}, we obtain that
\begin{equation}\label{intr1}
\frac1{p'}||Hf||_p\le ||H^*f||_p\le p||Hf||_p\quad\mbox{for} \quad 1<p<\infty
\end{equation}
(as usual, $p'=p/(p-1)).$

However, the constants in (\ref{intr1}) are not optimal. The
objective of this paper is to find optimal constants. Our main
result is the following theorem.

\begin{teo}\label{Main1} Let $f\in \mathcal M^+(\R_+)$ and let  $1<p<\infty.$  Then
\begin{equation}\label{main1}
(p-1)||Hf||_p \le ||H^*f||_p\le (p-1)^{1/p}||Hf||_p
\end{equation}
if $1<p\le 2$, and
\begin{equation}\label{main2}
(p-1)^{1/p}||Hf||_p \le ||H^*f||_p\le (p-1)||Hf||_p
\end{equation}
if $2\le p<\infty.$ All constants in (\ref{main1}) and (\ref{main2})
are the best possible.
\end{teo}

Clearly, the problem on relations between various norms of Hardy
operator and its dual is of independent interest (cf. \cite{CGMP2}).
At the same time, this problem has an equivalent formulation in
terms of the difference operator $H\f-\f.$

Let $\f$ be a nonincreasing and
nonnegative function on $\R_+$ such that $\f(+\infty)=0.$ The
quantity $H\f-\f$
plays an important role in Analysis (see \cite{BS}, \cite{CGMP1},
\cite{CGMP2}, \cite{K2007}, \cite{KS} and references therein). It is
well known that the norms $||H\f-\f||_p$ and $||\f||_p$
 ($1<p<\infty$) are equivalent
(see \cite[p. 384]{BSh}). However, the {\it sharp} constant is known only in the
following inequality.

Let $\f$ be a nonincreasing and
nonnegative function on $\R_+$.  Then for any
$p\ge 2$
\begin{equation}\label{intr2}
||H\f-\f||_p\le (p-1)^{-1/p}||\f||_p,
\end{equation}
and the constant is optimal.

This result  was obtained in \cite{KS} for $p=2k\,\, (k\in \N)$ and
in \cite{BS} for all $p\ge 2$ (we observe that (\ref{intr2}) is a special case of the inequality proved in
\cite{BS} for weighted $L^p-$norms).

 We shall show that inequality
(\ref{intr2}) is equivalent to the first inequality in (\ref{main2}):
\begin{equation}\label{intr3}
||Hf||_p \le (p-1)^{-1/p}||H^*f||_p, \quad 2\le p<\infty.
\end{equation}
Thus, (\ref{intr3}) can be derived from (\ref{intr2}).  However,
below we give a simple direct proof of (\ref{intr3}). Moreover,
Theorem \ref{Main1} has the following equivalent form.
\begin{teo}\label{Main2} Let $\f$ be a nonincreasing and
nonnegative function on $\R_+$ such that $\f(+\infty)=0$ and let $1<p<\infty.$ Then
\begin{equation}\label{main30}
(p-1)||H\f-\f||_p \le ||\f||_p\le (p-1)^{1/p}||H\f-\f||_p
\end{equation}
if $1<p\le 2$, and
\begin{equation}\label{main40}
(p-1)^{1/p}||H\f-\f||_p \le ||\f||_p\le (p-1)||H\f-\f||_p
\end{equation}
if $2\le p<\infty.$ All constants in (\ref{main30}) and (\ref{main40})
are the best possible.
\end{teo}

\section{Proofs of  main results}

{\it Proof of Theorem \ref{Main1}.} Taking into account
(\ref{intr1}), we may assume that $Hf$ and $H^*f$ belong to
$L^p(\R_+).$ We may also assume that $f(x)>0$ for all $x\in \R_+.$
Denote
\begin{equation*}
I_p=\int_0^\infty\left(\frac1x\int_0^x f(t)\,dt\right)^p\,dx.
\end{equation*}
Since $Hf\in L^p(\R_+),$ we have
$$
Hf(x)=o(x^{-1/p})\quad\mbox{as}\quad x\to 0+\quad\mbox{or}\quad x\to +\infty.
$$
Thus, integrating by parts, we obtain
\begin{equation}\label{main5}
I_p=p'\int_0^\infty x^{1-p}f(x)\left(\int_0^x f(t)\,dt\right)^{p-1}\,dx.
\end{equation}

Further, set
\begin{equation}\label{main6}
I_p^*=\int_0^\infty\left(\int_t^\infty \frac{f(x)}{x}\,dx\right)^p\,dt.
\end{equation}
First we shall prove that
\begin{equation}\label{main7}
(p-1)I_p\le I_p^*\quad\mbox{if}\quad 2\le p<\infty
\end{equation}
and
\begin{equation}\label{main8}
I_p^*\le (p-1)I_p \quad\mbox{if}\quad 1< p\le 2.
\end{equation}
Set
$$
\Phi(t,x)=\int_t^x\frac{f(u)}{u}\,du, \, 0<t\le x,
$$
and $G(t,x)= \Phi(t,x)^p.$ Since $G(t,t)=0,$ we have
$$
\left(\int_t^\infty \frac{f(x)}{x}\,dx\right)^p = \int_t^\infty G_x'(t,x)\,dx
=p \int_t^\infty \frac{f(x)}{x} \Phi(t,x)^{p-1}\,dx.
$$
Thus, by Fubini's theorem,
$$
I_p^*=p\int_0^\infty \int_t^\infty \frac{f(x)}{x} \Phi(t,x)^{p-1}\,dx\,dt
$$
\begin{equation}\label{main9}
=p\int_0^\infty \frac{f(x)}{x}\int_0^x \Phi(t,x)^{p-1}\,dtdx.
\end{equation}
On the other hand, Fubini's theorem gives that
$$
\int_0^x f(t)\,dt=\int_0^x \Phi(t,x)\,dt.
$$
Hence, by (\ref{main5}),
\begin{equation}\label{main10}
I_p=p'\int_0^\infty x^{1-p}f(x)\left(\int_0^x \Phi(t,x)\,dt\right)^{p-1}\,dx.
\end{equation}

Comparing (\ref{main5})
with (\ref{main6}), we see that $I_2=I_2^*.$ In what follows we assume that $p\not=2.$

Let $p>2.$ Then by H\"older's inequality
$$
\left(\int_0^x \Phi(t,x)\,dt\right)^{p-1}\le x^{p-2}\int_0^x \Phi(t,x)^{p-1}\,dt.
$$
Thus, by (\ref{main9}) and (\ref{main10}),
$$
I_p\le p'\int_0^\infty\frac{f(x)}{x}\int_0^x \Phi(t,x)^{p-1}\,dtdx=\frac{I^*_p}{p-1},
$$
and we obtain (\ref{main7}).

Let now $1<p< 2.$ Applying H\"older's inequality, we get
$$
\int_0^x \Phi(t,x)^{p-1}\,dt\le x^{2-p}\left(\int_0^x \Phi(t,x)\,dt\right)^{p-1}.
$$
Thus, by (\ref{main9}) and (\ref{main10}),
$$
I_p^*\le p\int_0^\infty x^{1-p}f(x)\left(\int_0^x \Phi(t,x)\,dt\right)^{p-1}\,dx=(p-1)I_p,
$$
and we obtain (\ref{main8}).

Inequalities (\ref{main7}) and (\ref{main8}) imply the first
inequality in (\ref{main2}) and the second inequality in
(\ref{main1}), respectively.

Now we shall show that
\begin{equation}\label{main11}
I_p^*\le (p-1)^pI_p\quad\mbox{if}\quad 2< p<\infty
\end{equation}
and
\begin{equation}\label{main12}
(p-1)^pI_p\le I_p^* \quad\mbox{if}\quad 1< p<2.
\end{equation}

Observe that by our assumption ($f>0$ and $H^*f\in L^p(\R^*))$,
$$
0<\int_t^\infty \frac{f(x)}{x}\,dx <\infty \quad\mbox{for all} \quad t>0.
$$
Thus, for any $q>0$ we have
\begin{equation}\label{main13}
\left(\int_t^\infty \frac{f(x)}{x}\,dx\right)^q=q\int_t^\infty \frac{f(x)}{x}\left(\int_x^\infty \frac{f(u)}{u}\,du\right)^{q-1}\,dx.
\end{equation}
Applying  this equality with $q=p$ in (\ref{main6}) and using Fubini's theorem, we
obtain
\begin{equation}\label{main14}
I_p^*=p\int_0^\infty f(x)\left(\int_x^\infty \frac{f(u)}{u}\,du\right)^{p-1}\,dx.
\end{equation}
Further,  apply (\ref{main13}) for $q=p-1$ and use again Fubini's
theorem. This gives
$$
I_p^*=p(p-1)\int_0^\infty f(x)\int_x^\infty \frac{f(u)}{u}\left(\int_u^\infty \frac{f(v)}{v}\,dv\right)^{p-2}\,du\,dx
$$
$$
= p(p-1)\int_0^\infty\frac{f(u)}{u}\left(\int_u^\infty \frac{f(v)}{v}\,dv\right)^{p-2}\int_0^u f(x)\,dx\,du.
$$
 Set
$$
\f(u)=\frac{f(u)^{1/(p-1)}}{u}\int_0^u f(x)\,dx
$$
and
$$
\psi(u)=f(u)^{(p-2)/(p-1)}\left(\int_u^\infty \frac{f(x)}{x}\,dx\right)^{p-2}
$$
(recall that $f>0$). Then we have
\begin{equation}\label{main15}
I_p^*=p(p-1)\int_0^\infty \f(u)\psi(u)\,du.
\end{equation}
Furthermore, by (\ref{main5}),
\begin{equation}\label{main16}
\int_0^\infty \f(u)^{p-1}\,du=\int_0^\infty \frac{f(u)}{u^{p-1}}\left(\int_0^u f(x)\,dx\right)^{p-1}\,du= \frac{I_p}{p'},
\end{equation}
and by (\ref{main14}),
\begin{equation}\label{main17}
\int_0^\infty\psi(u)^{(p-1)/(p-2)}\,du= \int_0^\infty f(u)\left(\int_u^\infty \frac{f(x)}{x}\,dx\right)^{p-1}\,du=\frac{I_p^*}{p}
\end{equation}
for any $p>1, \, p\not=2.$

 Let $p>2.$ Applying in (\ref{main15}) H\"older's inequality with the
exponent $p-1$ and taking into account equalities (\ref{main16}) and
(\ref{main17}), we obtain
$$
I_p^*\le p(p-1)\left(\frac{I_p}{p'}\right)^{1/(p-1)}\left(\frac{I_p^*}{p}\right)^{(p-2)/(p-1)}.
$$
This implies   (\ref{main11}),  which is the second   inequality in
(\ref{main2}).

Let now $1< p<2.$ Applying in (\ref{main15}) H\"older's inequality with the exponent
$p-1\in (0,1)$ (see \cite[p. 140]{HLP}), and using equalities
(\ref{main16}) and (\ref{main17}), we get
$$
I_p^*\ge p(p-1)\left(\frac{I_p}{p'}\right)^{1/(p-1)}\left(\frac{I_p^*}{p}\right)^{(p-2)/(p-1)}.
$$
Thus,
$$
(I_p^*)^{1/(p-1)}\ge(p-1)^{p/(p-1)}I_p^{1/(p-1)}.
$$
This implies   (\ref{main12}),  which is the first  inequality in
(\ref{main1}).

It remains to show that the constants in (\ref{main1}) and
(\ref{main2}) are optimal. First, set
$f_\e(x)=\chi_{[1,1+\e]}(x)\,\, (\e>0).$ Then
$$
||Hf_\e||_p^p
=\int_1^{1+\e} x^{-p}(x-1)^p\,dx + \e^p\int_{1+\e}^\infty x^{-p}\,dx.
$$
Thus,
$$
\frac{\e^p(1+\e)^{1-p}}{p-1}\le ||Hf_\e||_p^p\le \frac{\e^p(1+\e)^{1-p}}{p-1} +\e^{p+1}.
$$
Further,
$$
\begin{aligned}
||H^*f_\e||_p^p &= \int_0^1\left(\int_1^{1+\e}\frac{dt}{t}\right)^p\,dx
+\int_1^{1+\e}\left(\int_x^{1+\e}\frac{dt}{t}\right)^p\,dx\\&= (\ln (1+\e))^p+\int_1^{1+\e}\left(\ln\frac{1+\e}{x}\right)^p\,dx.
\end{aligned}
$$
Thus,
$$
(\ln (1+\e))^p\le ||H^*f_\e||_p^p\le (\ln (1+\e))^p(1 +\e).
$$
Using these estimates, we obtain that
$$
\lim_{\e\to 0+}\frac{||Hf_\e||_p}{||H^*f_\e||_p}=(p-1)^{-1/p}.
$$
It follows that the constants in the right-hand side of (\ref{main1}) and the left-hand side of (\ref{main2}) cannot be improved.

Let $1<p<2.$ Set $f_\e(x)=x^{\e-1/p}\chi_{[0,1]}(x)\,\, (0<\e<1/p).$
Then
$$
||Hf_\e||_p^p\ge \int_0^1\left(\frac1x\int_0^x t^{\e-1/p}\,dt\right)^p\,dx=\frac{p^p}{\e p(p-1+\e p)^p}.
$$
On the other hand,
$$
||H^*f_\e||_p^p\le \left(\frac1p-\e\right)^{-p}\int_0^1 x^{(\e-1/p)p}\,dx=\frac{p^p}{\e p(1-\e p)^p}.
$$
Hence,
$$
\varliminf_{\e\to 0+}\frac{||Hf_\e||_p}{||H^*f_\e||_p}\ge \frac1{p-1}.
$$
This implies that the constant in the left-hand side of
(\ref{main1}) is optimal.

Let now $p>2.$ Set $f_\e(x)=x^{-\e-1/p}\chi_{[1,+\infty)}(x)\,\,
(0<\e<1/p').$ Then
$$
||H^*f_\e||_p^p\ge\int_1^\infty\left(\int_x^\infty\frac{dt}{t^{1+1/p+\e}}\right)^p\,dx=\frac{p^p}{\e p(1+\e p)^p}
$$
and
$$
||Hf_\e||_p^p\le \int_1^\infty\left(\frac1x\int_0^x\frac{dt}{t^{1/p+\e}}\right)^p\,dx=\frac{p^p}{\e p(p-1-\e p)^p}.
$$
Thus,
$$
\varliminf_{\e\to 0+}\frac{||H^*f_\e||_p}{||Hf_\e||_p}\ge p-1.
$$
This shows that the constant in the right-hand side of (\ref{main2})
is the best possible. The proof is completed.

\vskip 8pt

\begin{rem} We emphasize that in Theorem \ref{Main1} we do not assume that $f$ belongs to $L^p(\R_+).$ It is clear
that the condition $Hf\in L^p(\R_+)$ does not imply
that $f\in L^p(\R_+).$ For example, let $f(x)=|x-1|^{-1/p}\chi_{[1,2]}(x),$ $p>1.$ Then
$$
Hf(x)=0\quad\mbox{for}\quad x\in [0,1]\quad\mbox{and}\quad Hf(x)\le \frac{p'}{x}\quad\mbox{for}\quad x\ge 1.
$$
Thus, $Hf\in L^p(\R_+)$, but $f\not\in L^p(\R_+).$
\end{rem}

Now we shall  show that Theorems \ref{Main1} and \ref{Main2} are
equivalent. First we observe that without loss of generality we may
assume that a function $\f$ in Theorem \ref{Main2} is locally
absolutely continuous on $\R_+.$ Indeed, let $\f$ be a nonincreasing
and nonnegative function on $\R_+$ such that $\f(+\infty)=0.$ Set
$$
\f_n(x)=n\int_x^{x+1/n} \f(t)\,dt \quad (n\in \N).
$$
Then functions $\f_n$ are nonincreasing, nonnegative, and locally
absolutely continuous on $\R_+.$ Besides, the sequence $\{\f_n(x)\}$
increases for any $x\in \R_+$ and converges to $\f(x)$ at every
point of continuity of $\f.$ By the monotone convergence theorem,
$H\f_n(x)\to H\f(x)$ as $ n\to \infty$ for any $x\in \R_+$, and
$||\f_n||_p\to ||\f||_p.$ Furthermore, in Theorem \ref{Main2} we may
assume that  $\f\in L^p(\R_+)$ (in conditions of this theorem the
norms  $||H\f-\f||_p$ and $||\f||_p$ are equivalent \cite[p.
384]{BSh}). Using this assumption, Hardy's inequality,  and
 the dominated convergence theorem, we obtain that $||H\f_n-\f_n||_p\to
||H\f-\f||_p$.

Let $\f$ be a nonincreasing,
nonnegative, and locally absolutely continuous function  on $\R_+$ such that $\f(+\infty)=0.$  Then
$$
H\f(x)-\f(x)=\frac1x\int_0^x[\f(t)-\f(x)]\,dt
$$
$$
=\frac1x\int_0^x\int_t^x|\f'(u)|\,dudt=\frac1x\int_0^x u|\f'(u)|\,du.
$$
Set $u|\f'(u)|=f(u).$ Since $\f(+\infty)=0,$ we have
$$
\f(x)=\int_x^\infty |\f'(u)|\,du=\int_x^\infty \frac{f(u)}{u}\,du.
$$
Thus,
\begin{equation}\label{monot1}
H\f(x)-\f(x)=\frac1x\int_0^x f(u)\,du = Hf(x)
\end{equation}
and
\begin{equation}\label{monot2}
\f(x)=\int_x^\infty \frac{f(u)}{u}\,du=H^*f(x).
\end{equation}

Conversely, if $f\in \mathcal M^+(\R_+)$ and
$$
\int_0^x f(u)\,du<\infty \quad\mbox{for any}\quad x>0,
$$
we define $\f$ by (\ref{monot2}) and then we have equality
(\ref{monot1}). These arguments show the equivalence of Theorems
\ref{Main1} and \ref{Main2}.

\end{document}